\documentclass{ifacconf}

\usepackage{algorithm}
\usepackage{algpseudocode}

\usepackage[all]{xy}
\usepackage{graphicx}      
\usepackage{natbib}        
\usepackage{amsmath,amssymb,amsfonts,mathrsfs,stmaryrd}
\usepackage{array,multirow,makecell}
\usepackage[latin1]{inputenc}
\usepackage{tikz}

\usepackage{color}
\usepackage{bbm}

 \newtheorem{theorem}{Theorem}

\newtheorem{lemma*}{Lemma}
{}
{}

\newtheorem{remark}[theorem]{Remark}

\newcommand {\R}   {\mathbbm R}

\newcommand {\C}   {\mathbbm C}

\newcommand {\Ubar}   {\overline{\mathbb U}}

\newcommand {\N}   {\mathbbm N}
\newcommand {\Q}   {\mathbbm Q}

\newcommand{\makeremark}[2]{
  \newcommand{#1}[1]
    {
    \color{blue}
     $\longrightarrow$ \textsc{#2: }
     ##1
     $\longleftarrow$
    \color{black}
    }
}    

\makeremark{\FR}{Fabrice says}
\makeremark{\AQ}{Alban says}
\makeremark{\YB}{Yacine says}

\definecolor{1ST}{rgb}{1,0,0}
\definecolor{2ND}{rgb}{1,0.5,0}
\definecolor{3RD}{rgb}{1,0,1}

\newcommand{\shutup}[1]{}

\renewcommand{\leq}{\leqslant}  
\renewcommand{\geq}{\geqslant}

\makeatletter
\def\cramped                           
 {\parskip0pt\@topsep0pt       
  \itemsep0pt\parsep0pt
}        
\makeatother

\usepackage{stmaryrd}

\begin{document}
	\begin{frontmatter}
		
		\title{Computing effectively stabilizing controllers for a class of $n$D systems} 
		
		\thanks[footnoteinfo]{All but the third author were supported by the ANR-13-BS03-0005 (MS-DOS).}
		
		\author[First]{Yacine Bouzidi} 
		\author[Second]{Thomas Cluzeau}
		\author[Third]{Guillaume Moroz}
		\author[First]{Alban Quadrat}

		\address[First]{INRIA Lille-Nord Europe, Lille, France \\(e-mail: \{yacine.bouzidi,alban.quadrat\}@inria.fr).}
		\address[Second]{University of Limoges ; CNRS ; XLIM UMR 7252, Limoges, France (e-mail: thomas.cluzeau@unilim.fr).}
		\address[Third]{INRIA Nancy-Grand Est, Nancy, France \\(e-mail: guillaume.moroz@inria.fr)}

		\begin{abstract}
			
			In this paper, we study the internal stabilizability and internal stabilization problems for multidimensional ($n$D) systems. Within the fractional representation approach, a multidimensional system can be studied by means of matrices with entries in the integral domain of structurally stable rational fractions, namely the ring of rational functions which have no poles in the closed unit polydisc  $\Ubar^n= \left \{z=(z_1,\ldots,z_n) \in \C^n \ | \ |z_1|\leq 1,\ldots,|z_n|\leq 1 \right \}.$ 
			It is known that the internal stabilizability of a multidimensional system can be investigated by studying a certain polynomial ideal $I=\langle p_1,\ldots,p_r \rangle$ that can be explicitly described in terms of the transfer matrix of the plant. More precisely the system is stabilizable if and only if $V(I)=\{z \in \C^n \; | \; p_1(z)=\cdots=p_r(z)=0 \} \cap \Ubar^n = \emptyset$. In the present article, we consider the specific class of linear $n$D systems (which includes the class of 2D systems) for which the ideal $I$ is zero-dimensional, i.e., the $p_i$'s have only a finite number of common complex zeros. We propose effective symbolic-numeric algorithms for testing if $V(I) \cap \Ubar^n = \emptyset$, as well as for computing, if it exists, a stable polynomial $p \in I$ which allows the effective computation of a stabilizing controller. We illustrate our algorithms through an example and finally provide running times of prototype implementations for 2D and 3D systems. 
		\end{abstract}
		
		
	\end{frontmatter}
	
	\section{Introduction}\label{sec:intro}
	
	Multidimensional or $n$D systems (\cite{bose84}) are systems of functional equations whose unknown functions depend on $n$ independent variables. The stabilizability and stabilization problems are fundamental issues in the study of multidimensional  systems in control theory. Nowadays, the problem is well-understood in the case of 1D systems whereas progress for $n$D systems with $n\geq2$ are rather slow. One approach for handling stabilizability or stabilization issues in systems theory is the {\em fractional representation approach} (\cite{vidyasagar2011control}) in which a plant is represented by its transfer matrix $P \in K^{q \times r}$ where $K=\R(z_1, \ldots, z_n)$. This transfer matrix admits a left factorization $P=D^{-1}\,N$ (also called fractional representation of $P$), where the matrices $D \in A^{q \times q}$ satisfying $\det(D) \neq 0$ and $N \in A^{q \times r}$  have entries in the integral domain $A=\R(z_1, \ldots, z_n)_S$ of structurally stable rational fractions, namely the ring of rational functions in $z_1, \ldots, z_n$ which have no poles in the closed unit polydisc of $\C^n$ defined by: $$\Ubar^n= \left \{z=(z_1,\ldots,z_n) \in \C^n \ | \ |z_1|\leq 1,\ldots,|z_n|\leq 1 \right \}.$$
	Introducing the matrix $R = (D \quad -N) \in A^{q \times (q+r)}$, it is known (see \cite{quadrat2003fractional1,quadrat2003fractional2}) that the multidimensional system given by the transfer matrix $P=D^{-1}\,N$ is then internally stabilizable if and only if the $A$-module $A^{1 \times (q+r)}/\overline{A^{1 \times q}\,R}$ is a projective $A$-module of rank $r$, where the closure $\overline{A^{1 \times q}\,R}$ of $A^{1 \times q}\,R$ in $A^{1 \times (q+r)}$ is defined by: $\overline{A^{1 \times q}\,R} = \{\lambda \in A^{1 \times (q+r)} \; | \; \exists \ a \in A \setminus \{0\}: a\,\lambda \in  A^{1 \times q}\,R \}.$
	This projectivity condition is in turn equivalent to the fact that the reduced minors of the matrix $R$ do not have common zeros in $\Ubar^n$ (see also \cite{lin1998}). In other terms, if we denote by $p_1, \ldots, p_r$ the \emph{reduced minors} of $R$, i.e., the $q\times q$ minors of $R$ divided by their gcd, by $I=\langle p_1,\ldots,p_r \rangle$ the polynomial ideal generated by the $p_i$'s, and by $V(I)=\{z \in \C^n \; | \; p_1(z)=\cdots=p_r(z)=0\}$ the associated algebraic variety, then the system is internally stabilizable if and only if $V(I) \cap  \Ubar^n = \emptyset$.
	
	The first contribution of the present paper is to provide an effective algorithm for testing the stabilizability condition $V(I) \cap  \Ubar^n = \emptyset$ for the class of $n$D systems for which the ideal $I$ is zero-dimensional, i.e., the $p_i$'s have only a finite number of common complex zeros (i.e., $V(I)$ consists of a finite number of complex points). Note that this class includes the class of 2D systems. Our main idea is to take advantage of the univariate representation for zero-dimensional ideals (\cite{canny1988some,becker1996radical,alonso1996zeros,rouillier1999solving}). This concept, which can be traced back to \cite{kronecker1882grundzuge}, yields a one-to-one correspondance between the elements of $V(I)$ and the zeros of a univariate polynomial $f$. Numerical techniques can thus be applied to compute certified numerical approximations of the roots of $f$ and then of those of $V(I)$.
	
	In the case of a stabilizable plant, the next step consists in computing a stabilizing controller which can be achieved by computing a stable (i.e., devoid from zeros in $\Ubar^n$) polynomial $s \in I$ (see \cite{lin19882D}). The {\em polydisc Nullstellensatz}, proved by Bridges, Mines, Richman and Schuster (see \cite{bridges2004}), shows that the existence of a stable polynomial $s \in I$ is equivalent to $V(I)\cap \Ubar^n = \emptyset$. Several proofs of this result have been investigated in the literature, mainly for the case where $I$ is a zero-dimensional ideal (see \cite{raman1986constructive,lin19882D,bisiacco1986controller,guiver1995causal} and \cite{xu1994} for instance). Nevertheless none of them is effective in the sense that it provides an algorithm for computing $s$ using calculations that can be performed in an exact way by a computer. Indeed, starting from a set of polynomials with rational coefficients ($I \subset \Q[z_1,\ldots,z_n]$),  these algorithms are built on \emph{spectral factorization}, i.e., factorization of polynomials in  $\Q[z]$ into stable and instable factors. For irreducible polynomials in $\Q[z]$, this factorization requires the explicit computation of the complex roots of the polynomials, which can be done only approximately. This leads to approximate (stable) polynomials that do not belong to the polynomial ideal. As a consequence, these algorithms are able to solve the aforementioned problem only for few simple systems (see Sections~\ref{sec:stabilization} and \ref{sec:experiments} for details).
	
	
	Our second contribution is to provide an effective algorithm for computing a stable polynomial $s=\sum_{i=1}^r u_i \, p_i \in I$ for the class of systems for which $I \subset \Q[z_1,\ldots,z_n]$ is a zero-dimensional ideal. Our symbolic-numeric method roughly follows the lines of that proposed in \cite{xu1994} but once again we take advantage of the univariate representation of zero-dimensional ideals (\cite{rouillier1999solving}) to control the numeric precision required to achieve our goal. 
	
	The paper is organized as follows. In Section~\ref{sec:solving}, we recall some classical  computer algebra results on the complex zeros of polynomials and polynomial systems. We also introduce the univariate representation of zero-dimensional ideals which will be our main tool in what follows. In Section~\ref{sec:stabilizability}, we provide an effective stabilizability test, i.e., an algorithm for testing whether a zero-dimensional ideal intersects the closed unit polydisc. In Section~\ref{sec:stabilization}, we provide an effective \emph{polydisc Nullstellensatz} namely a symbolic-numeric method for computing, if it exists, a stable polynomial in a zero-dimensional polynomial ideal. Finally, in Section~\ref{sec:experiments}, we illustrate our methods on one example and show some running times of prototype implementations.

	\section{Preliminaries on algebraic systems}\label{sec:solving}
	
	In this section, we introduce some notations and we recall some classical material about the computation of certified numerical approximations of the complex zeros of polynomials and polynomial systems.
	
	The {\em bit-size} of an integer is the number of bits in its representation and for a rational number (resp., a polynomial with rational coefficients) the term {\em bit-size} refers to 
	the maximum bit-size of its numerator and denominator (resp., of its coefficients). For a complex number $z \in \C$, we denote by $\Re(z) \in \R$ (resp., $\Im(z) \in \R$) its real (resp., imaginary) part. If $z_1,\,z_2 \in \C$, we write {\em $z_1 < z_2$} if both $\Re(z_1) < \Re(z_2)$ and $\Im(z_1) < \Im(z_2)$. For $z_1,\,z_2 \in \C$ such that $z_1 < z_2$, we shall consider the {\em axes-parallel open box} or {\em box} for short $B=B(z_1,z_2)=\{z \in \C \ | \ z_1 < z < z_2 \}$ and its {\em width} is defined by $w(B)=\max\{|\Re(z_2-z_1)|,|\Im(z_2-z_1)|\}$. We also introduce the non-negative real number  $|B|=\max\{|\Re(z_1)|,|\Im(z_1)|,|\Re(z_2)|,|\Im(z_2)|\}$. The box $B$ is said to be of {\em rational endpoints} if $z_1,$ $z_2$ have rational real and imaginary parts, i.e, $\Re(z_1),\,\Im(z_1),\,\Re(z_2),\,\Im(z_2) \in \Q$. Finally,  the box $B$ is called {\em isolating} for a given polynomial $f \in \Q[z]$ if it contains exactly one complex zero of $f$.  
	
	The following result concerns the isolation of the complex zeros of a univariate polynomial. We refer, for instance, to \cite{sagraloff2009} for more details.
	\begin{lem}\label{thm:isolation}
		Let $f \in \Q[z]$ be a squarefree polynomial of degree $d$. Then, for all $\epsilon>0$, one can compute disjoint axes-parallel open boxes $B_1,\ldots,B_d$, with rational endpoints such that each $B_i$ contains exactly one complex root of $f$ and satisfies $w(B_i)\leq \epsilon$.	
	\end{lem}
	
	In the algorithms given in Sections~\ref{sec:stabilizability} and \ref{sec:stabilization} below, we shall use a routine called {\tt Isolate} which takes as input a univariate polynomial $f$, a box $B$, and a precision $\epsilon >0$ and computes isolating boxes $B_1,\ldots,B_l$ with rational endpoints for the complex roots of $f$ that lie inside the given box $B$ and such that $\max_{i=1,\ldots,l} w(B_i) \leq \epsilon$. If $B$ (resp., $\epsilon$) is not specified in the input, we consider all complex roots in $\C$  (resp., the boxes are computed up to a sufficient precision for isolation). 
	
	Let us now recall a standard property about width expansion through interval arithmetic in polynomial evaluation. Here we consider {\em exact interval arithmetic}, that is, the arithmetic operations on the interval endpoints are considered exact (see \cite{alefeld2012introduction}). If $f \in \Q[x_1,x_2]$ is a bivariate polynomial of two real variables $x_1$ and $x_2$ and $B$ a box, we denote by $\square f(B)$ the interval that results from the evaluation of the polynomial $f$ at the box $B$ using interval arithmetic.
	\begin{lem}[\cite{cheng2010topology}, Lemma 8]\label{thm:interval}
		Let $B$ be a box with rational endpoints satisfying $|B| \leqslant 2^{\sigma}$ and let $f \in \Q[x_1,x_2 ]$ be a bivariate polynomial of two real variables $x_1$ and $x_2$ of degree $d$ with coefficients of bit-size $\tau$. Then, $f$ can be evaluated at the box $B$ by interval arithmetic into an interval $\square f(B)$ of width at most $2^{\tau + d\,\sigma+1} \,d^3\,w(B)$. 
	\end{lem}
	In particular, a direct consequence of Lemma~\ref{thm:interval} is that if $w(B) \leq \epsilon \, 2^{ -\tau - d\,\sigma-1- 3 \log_2(d)}$, then we have $w(\square f(B)) \leq \epsilon$.	 
	
	We now consider a set of polynomials $p_1,\ldots,p_r$ in $\Q[z_1,\ldots,z_n]$. We denote by $I=\langle p_1,\ldots,p_r\rangle$ the ideal generated by the $p_i$'s and by $$V(I)=\{z \in \C^n \; | \; p_1(z)=\cdots=p_r(z)=0\} \subseteq \C^n,$$ the complex variety of their common zeros. In the sequel, we shall always assume that the ideal $I$ under consideration is a \emph{zero-dimensional} ideal, that is, that the $p_i$'s have only a finite number of common complex zeros, i.e., $V(I)$ consists of a finite number of complex points. Methods for computing certified numerical approximations for the elements of $V(I)$ usually proceed in two steps. First, a formal representation (as for instance, a \emph{Gröbner basis}, a \emph{triangular decomposition} or a \emph{univariate representation}) of the set $V(I)$ is computed. Then, this formal representation is used to compute, more or less easily, numerical approximations of the elements of $V(I)$. The convenient representation of the variety of zero-dimensional ideals that we shall use in the sequel is the so-called {\em univariate representation} introduced by Rouillier in \cite{rouillier1999solving}.  
	\begin{defn}\label{def:UR}
		With the previous notation and assumptions, a \emph{univariate representation} of $V(I)$ is the datum of a linear form $t = a_1 \, z_1 + \cdots + a_n \,z_n$, with $a_1,\ldots,a_n \in \Q$ as well as $n+1$ univariate polynomials $f,\, g_{z_1}, \ldots, g_{z_n} \in \Q[t]$  such that the following two applications 
		$$\begin{array}{ccl}
		V(I) & \longrightarrow & V(f)=  \{z \in \C \; | \; f(z)=0\}\\
		z=(z_1,\ldots,z_n) &  \longmapsto & a_1 \, z_1 + \cdots + a_n \,z_n,
		\end{array}$$ 
		and
		$$
		\begin{array}{ccl}
		V(f) & \longrightarrow & V(I)\\
		z &  \longmapsto & (g_{z_1}(z),\ldots,g_{z_n}(z)),
		\end{array} $$
		provide a one-to-one correspondence between the elements of $V(I)$ and the zeros of $f$.
	\end{defn}	
	The univariate representation then defines a bijection between the zeros of $I$ and those of $f$ that preserves the multiplicities and the real zeros. One of its main advantage is that it permits to study the variety $V(I)$, i.e., the set of common zeros of the polynomials $p_i$'s of $n$ variables through the roots of the univariate polynomial $f$. In particular, if we are interested in computing certified numerical approximations of the complex elements of $V(I)$, a first step consists, in isolating the complex roots of the polynomial $f$ (using for instance the algorithm in \cite{sagraloff2009}). Then, numerical approximations of the coordinates $z_1,\ldots,z_n$ of $z \in V(I)$ can be obtained by evaluating $g_{z_1},\ldots,g_{z_n}$ at the obtained numerical boxes. Lemmas~\ref{thm:isolation} and \ref{thm:interval} above show that we can thus obtain a given precision for the elements of $V(I)$.
	
	The computation of a univariate representation of $V(I)$ can be done by pre-computing a Gröbner basis $G$ of the ideal $I$ \cite{cox1992ideals}, and then performing linear algebra calculations in the quotient vector space $\Q[z_1,\ldots,z_n]/I$ which admits a basis composed of the monomials that are irreducible modulo the Gröbner basis $G$. For more details, see \cite{rouillier1999solving}. For the specific case of bivariate algebraic systems, it should be stressed that very practically efficient algorithms exist for computing such a representation: see \cite{bouzidi2014resolution}. To achieve efficiency, these algorithms replace the ``costly" computation of a Gröbner basis of $I$ by \emph{resultants} and \emph{subresultants} computations, see \cite{basu2006algorithms}.
	
	In the algorithms given in Sections~\ref{sec:stabilizability} and \ref{sec:stabilization} below, we shall use a routine called {\tt Univ\_R} which takes as input the polynomials $p_1,\ldots,p_r$ and computes a univariate representation of the variety $V(I)$ of the zero-dimensional ideal $I=\langle p_1,\ldots,p_r\rangle$.
	
	\section{An effective stabilizability test}\label{sec:stabilizability}
	
	Using the fractional representation approach to multidimensional systems, in Section~\ref{sec:intro}, we have seen that a plant is stabilizable if and only if $V(I) \cap \Ubar^n = \emptyset$ for a certain polynomial ideal $I= \langle p_1,\ldots,p_r\rangle$ that can be explicitly described in terms of the transfer matrix of the plant. 
	
	Given polynomials $p_1,\ldots,p_r \in \Q[z_1,\ldots,z_n]$, the purpose of this section is to provide an effective algorithm to decide whether or not $V(I) \cap \Ubar^n = \emptyset$ where $I= \langle p_1,\ldots,p_r\rangle$ is a zero-dimensional ideal. To achieve this, we shall use a symbolic-numeric approach. We start by computing a univariate representation of $V(I)$. As explained in Section~\ref{sec:solving}, such a representation allows to describe formally the elements $z = (z_1,\ldots,z_n)$ of $V(I)$ as
	
	\begin{equation}\label{eq:rur}
	\{f(t) = 0, \quad z_1 = g_{z_1}(t), \quad \ldots, \quad z_n = g_{z_n}(t)\}
	\end{equation}
	
	where $f,g_{z_1},\ldots,g_{z_n} \in \Q[t]$. In what follows, the degree of $f$ is denoted by $d$, and those of $g_{z_1},\ldots,g_{z_n}$ are then smaller than $d$ (see \cite{rouillier1999solving}).
	
	Using a univariate representation, one can compute a set of hypercubes in $\R^{2n}$ isolating the elements of $V(I)$. Each coordinate is represented by a box in $\R^2$ obtained from the intervals containing its real and imaginary parts. Moreover, from Lemmas~\ref{thm:isolation} and \ref{thm:interval}, these hypercubes can be refined up to an arbitrary precision. We shall now consider the intersection between those hypercubes and the closed unit polydisc of $\C^n$ defined by: $$\Ubar^n= \left \{z=(z_1,\ldots,z_n) \in \C^n \ | \ |z_1|\leq 1,\ldots,|z_n|\leq 1 \right \}.$$

	Below, for any $g_{z_i}\in \Q[t]$, we shall denote by ${\cal C}(g_{z_i})$ the bivariate polynomial $\Re(g_{z_i})^2+\Im(g_{z_i})^2-1 \in \Q[x_1,x_2]$, where $\Re(g_{z_i})$ (resp., $\Im(g_{z_i})$) is the real (resp., complex) part of the polynomial resulting from $g_{z_i}(x_1+i\,x_2)$.

	From the definition of $\Ubar^n$, one can see that the situation is easier when $V(I)$ does not contain elements $z \in \C^n$ with $|z_i|=1$ for some $i \in \{1,\ldots,r\}$. Indeed, we have:
	\begin{thm}\label{thm:refine}
		With the previous notations, let us consider $z=(z_1,\ldots,z_n) \in V(I)$ such that, for all $i \in \{1,\ldots,n\}$, $|z_i|\neq 1$. Let $B$ be an isolating box for the root of $f$ corresponding to $z$ in the univariate representation of $V(I)$.
		Then, there exists $\epsilon>0$ such that if $w(B) \leq \epsilon$, then, for all $k \in \{1,\ldots,n\}$, the interval $\square{\cal C}(g_{z_k})(B)$ does not contain zero. 
	\end{thm} 
	\begin{pf}
		Let $d$ (resp., $\tau$) denote an upper bound on the degree (resp., bit-size) of the polynomials $g_{z_k}$, $k=1,\ldots,n$. For $k \in \{1,\ldots,n\}$, the real  (resp., imaginary) part of $g_{z_k}(x_1+i\,x_2)$ is a bivariate polynomial in $x_1$ and $x_2$ of degree (resp., bit-size) bounded by $d$ (resp., $d+\tau$). Consequently, the bivariate polynomial $ {\cal C}(g_{z_k}) \in \Q[x_1,x_2]$ has degree and bit-size respectively bounded by $2\,d$ and $(d+\tau)\,d$ respectively. Now, for $z=(z_1,\ldots,z_n) \in V(I)$, let $m = \min_{k=1,\ldots,n}||z_k|-1| > 0$ and let $B$ be an isolating box for the root of $f$ corresponding to $z$ in the univariate representation of $V(I)$ (\ref{eq:rur}), and such that $|B|\leqslant 2^{\sigma}$. From Lemma~\ref{thm:interval}, if we refine $B$ so that $w(B) \leq \epsilon$, where $\epsilon = m \, 2^{-(d+\tau)\,d - 2\,d\,\sigma-1- 3 \log_2(2\,d)}$, then, we have that for all $k \in \{1,\ldots,n\}$, the interval $\square {\cal C}(g_{z_k})(B)$ satisfies $w(\square{\cal C}(g_{z_k})(B)) \leq m$ so that, by definition of $m$, it does not contain zero.  
	\end{pf}
	Therefore, if $V(I)$ does not contain elements $z \in \C^n$ with one coordinate $z_i$ in the unit circle, one can easily test the stabilizability condition $V(I) \cap \Ubar^n = \emptyset$. Indeed, with the previous notations, we isolate the roots of $f$ inside boxes $B_1,\ldots,B_d$. Then, for $i \in \{1,\ldots,d\}$, we refine $B_i$ until, for all $k \in \{1,\ldots,n\}$, the interval $\square{\cal C}(g_{z_k})(B_i)$ does not contain zero. If one of the intervals $\square{\cal C}(g_{z_k})(B_i)$ is included in $\R_+$, we proceed to the next box $B_i$, otherwise we have found an element in $V(I) \cap \Ubar^n$ so that the system is certainly not stabilizable. After having investigated all the boxes $B_i$, we can then conclude about the stabilizability of the system. 
	
	We shall now consider the case where $V(I)$ contains (at least) one element having some coordinates on the unit circle. In this case, we cannot proceed numerically as before since if $|z_k|=1$, then, using the above notations, we cannot fulfill the condition that the interval $\square{\cal C}(g_{z_k})(B)$ does not contain zero.
	To guarantee the termination of the algorithm, we shall then have to compute, for each variable $z_k$, the number of elements 
	$z=(z_1,\ldots,z_n) \in V(I)$ satisfying $|z_k|=1$. 
	\begin{lem}\label{lem:cardroot=1}
		Let $I \subset \Q[z_1,\ldots,z_n]$ be a zero-dimensional ideal and $V(I)$ the associated algebraic variety. Then, for all $k \in \{1,\ldots,n\}$, one can compute the non-negative integer $l_k= \sharp \{z=(z_1,\ldots,z_n) \in V(I) \ | \ |z_k|=1\}$.
	\end{lem} 
	\begin{pf}
		Let $\{f(t),z_1-g_{z_1}(t),\ldots,z_n-g_{z_n}(t)\}$ be a univariate representation of $V(I)$ and $k \in \{1,\ldots,n\}$. Computing the resultant of the polynomials $f$ and $z_k-g_{z_k}$ with respect to the variable $t$ we get a univariate polynomial that can be written $r_k=\prod_{\alpha \in V(I)}{(z_k-\alpha_k)^{\mu_{\alpha_k}}}$, where the multiplicity $\mu_{\alpha_k}$ corresponds to $\sharp \{z \in V(I) \ | \ z_k=\alpha_k\}$. Then, using the classical \emph{Bistritz} test (see \cite{bistritz2002zero}), one can compute the number of complex roots counted with multiplicity of $r_k$ that lie on the unit circle and obtain the non-negative integer $l_k$.
	\end{pf}
	
	\begin{remark}
		An alternative to the Bistritz test, which we use in practice, consists in applying to the polynomial $r_k$, the Möbius transform $z_k \rightarrow  \frac{z_k-i}{z_k+i}$, which maps the complex unit circle $\Ubar$ to the real line $\R \cup \{\infty\}$. The number of complex roots of $r_k$ on the unit circle is then given as the number of real roots of the gcd of two polynomials $\Re$ and $\Im$, where $\Re$ (resp., $\Im$) is the real (resp., complex) part of the numerator of the rational fraction $r_k(\frac{z_k-i}{z_k+i})$.
	\end{remark}

	Using Lemma~\ref{lem:cardroot=1}, we can test the stabilizability condition $V(I) \cap \Ubar^n = \emptyset$ as follows. We start with the variable $z_1$. We refine the isolating boxes $B_1,\ldots,B_d$ for the roots of $f$ until exactly $l_1$ intervals $\square{\cal C}(g_{z_1})(B_i)$ contain zero. We throw away the boxes $B_i$'s such that the interval $\square{\cal C}(g_{z_1})(B_i)$ is included in $\R_+$ and we proceed similarly with the next variable $z_2$. If at some point we have thrown away all the boxes $B_i$'s, then the system is stabilizable. Otherwise the boxes which remain at the end of the process lead to elements of $V(I) \cap \Ubar^n$ so that the system is not stabilizable. 
	
	We summarize our symbolic-numeric method for testing stabilizability in the following {\tt IsStabilizable} algorithm. 
	
	\begin{algorithm}[h]\label{algo:stabilizability} 
		\caption{{\tt IsStabilizable}}
		{\bf Input:} A set of $r$ polynomials $p_1,\ldots,p_r \subset \Q[z_1,\ldots,z_n]$.
		
		{\bf Output:} {\tt True} if $V(\langle p_1,\ldots,p_r \rangle) \cap \Ubar^n = \emptyset$, else {\tt False}.
		
		{\bf Begin}
		
		$\diamond$ $\{f,g_{z_1},\ldots,g_{z_n}\}:=$ {\tt Univ\_R}($\{p_1,\ldots,p_r\}$);
		
		$\diamond$ $\{B_1,\ldots,B_d\}:=$ {\tt Isolate}($f$);
		
		$\diamond$ $L_B:=\{B_1,\ldots,B_d\}$ and $\epsilon:=\min_{i=1,\ldots,d}{w(B_i)}$;
		
		{\bf For} $k$ from $1$ to $n$ {\bf do}
		
		\hspace{0.3cm}  $\diamond$ $l_k:= \sharp \{z \in V(I) \ | \ |z_k|=1\}$ (see Lemma~\ref{lem:cardroot=1});
		
		
		\hspace{0.3cm} {\bf While} $\sharp \{i \ | \ 0 \in \square{\cal C}(g_{z_k})(B_i) \} > l_k$ {\bf do}
		
		\hspace{0.6cm} $\diamond$ $\epsilon := \epsilon/2$;
		
		\hspace{0.6cm} $\diamond$ For $i=1,\ldots,d$, set $B_i:=${\tt Isolate}($f,B_i,\epsilon$);
		
		\hspace{0.3cm} {\bf End While}
		
		\hspace{0.3cm} $\diamond$ $L_B:=L_B \setminus \{B_i \ | \ \square{\cal C}(g_{z_k})(B_i) \subset \R_+\}$;
		
		\hspace{0.3cm} $\diamond$ If $L_B = \{\}$, then {\bf Return} {\tt True} {\bf End If};
		
		{\bf End For}
		
		{\bf Return} {\tt False}.
		
		{\bf End}
	\end{algorithm}
	
	\section{An effective stabilization algorithm}\label{sec:stabilization}
	
	When a multidimensional system is stabilizable, one is then interested in computing a stabilizing controller. Within the fractional representation approach, this problem reduces to the following task (see Section~\ref{sec:intro} or \cite{lin19882D}): given a polynomial ideal $I=\langle p_1,\ldots,p_r \rangle$ satisfying $V(I) \cap \Ubar^n= \emptyset$, compute a polynomial $s \in \Q[z_1,\ldots,z_n]$ such that $s$ belongs to the ideal $I$ and $s$ is a {\em stable polynomial}, i.e., $s$ is devoid from zeros in $\Ubar^n$.  
	
	For an ideal $I \subset \Q[z_1,\ldots,z_n]$, the {\em polydisc Nullstellensatz} (see \cite{bridges2004}) asserts that the existence of a stable polynomial $s \in I$ is equivalent to $V(I)\cap \Ubar^n = \emptyset$. 
	\begin{thm}[Polydisc Nullstellensatz]\label{thm:PN}
		Let us consider a polynomial ideal $I=\langle p_1,\ldots,p_r\rangle \subset \Q[z_1,\ldots,z_n]$ such that $V(I) \cap \Ubar^n = \emptyset$. Then, there exist $r+1$ polynomials $s,u_1,\ldots,u_r$ in $\Q[z_1,\ldots,z_n]$ such that: 
		$$s=\sum_{i=1}^r u_i\,p_i \text{ and } V(s) \cap \Ubar^n = \emptyset.$$
	\end{thm} 
	Several proofs of Theorem~\ref{thm:PN} have been invesigated in the literature. Nevertheless none of them is effective in the sense that it provides an algorithm for computing $s$ and the {\em cofactors} $u_i$'s using calculations that can be performed in an exact way by a computer. In \cite{xu1994}, the authors study $2$D systems (i.e., $n=2$) for which the ideal $I$ under consideration is zero-dimensional. The idea of their method for computing a stable polynomial 
	$s \in I$ is to compute univariate elimination polynomials $r_{z_1} \in \Q[z_1]$ and $r_{z_2} \in \Q[z_2]$ with respect to each variable $z_1$ and $z_2$ and to factorize them into a stable and an unstable factor, i.e., for $i=1, \,2$, $r_{z_i}=r_{z_i}^{(s)}\,r_{z_i}^{(u)}$, where the roots of $r_{z_i}^{(s)}$ (resp., $r_{z_i}^{(u)}$) are outside (resp., inside) the closed unit disc $\Ubar$, and then, to compute the stable polynomial as $s=r_{z_1}^{(s)}\,r_{z_2}^{(s)}$.
	However, this approach presents a major drawback with respect to the effectiveness aspect. Indeed, when the elimination polynomial $r_{z_1}$ (resp., $r_{z_2}$)  is an irreducible polynomial in $\Q[z_1]$ (resp. $\Q[z_2]$), its stable factor $r_{z_1}^{(s)}$ (resp., $r_{z_2}^{(s)}$) could not be computed exactly since it will have coefficients in $\C$, and thus,  only an approximation of this polynomial can be obtained. As a consequence, the polynomial $s=r_{z_1}^{(s)}\,r_{z_2}^{(s)}$ will not belong to the ideal $I$.  	

	In the sequel, we present a symbolic-numeric algorithm for computing $s$ and the $u_i$'s that follows roughly the approach of \cite{xu1994} while we provide a way for tackling the effectiveness issue. Our main ingredient is the univariate representation of zero-dimensional ideals which allows us to compute and refine approximate factorizations over $\Q$ of the elimination polynomials $r_{z_i}$'s.
	
	Let $I=\langle p_1,\ldots,p_r \rangle \subset \Q[z_1,\ldots,z_n]$ be a zero-dimensional ideal such that $V(I) \cap \Ubar^n= \emptyset$. For simplicity reasons, in what follows, we further assume that the ideal $I$ is \emph{radical}, i.e., $\sqrt{I}=\{p \in \Q[z_1,\ldots,z_n] \ | \ \exists m \in \N^*, \, p^m \in I\}=I$. The elements $z = (z_1,\ldots,z_n)$ of $V(I)$ are given by a univariate representation
	$$ \{f(t) = 0, \quad z_1 = g_{z_1}(t), \quad \ldots, \quad z_n = g_{z_n}(t)\},$$
	where $t=a_1\,z_1+\cdots+a_n\,z_n$, $a_k\in \Q$ for $k=1,\ldots,n$, $f,g_{z_1},\ldots,g_{z_n} \in \Q[t]$, and $\deg(f) = d$. Since $I$ is a radical ideal, the polynomial $f$ is a squarefree polynomial and $f(t)=\prod_{i=1}^d{(t-\gamma_i)}$ for distincts $\gamma_1\ldots,\gamma_d \in \C$.  Moreover, from  Definition~\ref{def:UR}, if we introduce the polynomial ideal $I_r=\langle f(t), z_1- g_{z_1}(t),\ldots,z_n- g_{z_n}(t) \rangle \subset \Q[t,z_1,\ldots,z_n]$, then we have $I_r = I \cap \langle t-\sum_{k=1}^{n}{a_k\,z_k} \rangle$. In particular, if $p(t,z_1,\ldots,z_n) \in I_r$, then $p(\sum_{k=1}^{n}{a_k\,z_k},z_1,\ldots,z_n) \in I$.
	
	Let us first explain how we can compute approximations $\tilde{r}_{z_k}^{(s)}$ of the stable polynomials $r_{z_k}^{(s)}$ appearing in the method of \cite{xu1994} sketched above. Using Lemma~\ref{thm:isolation}, we can compute a set of boxes  $B_1,\ldots,B_d$ with rational endpoints, isolating the distinct complex roots $\gamma_1,\ldots,\gamma_d$ of $f$. Then, according to the stabilizability condition $V(I) \cap \Ubar^n= \emptyset$, for all $i \in \{1,\ldots,d\}$, the box $B_i$ can be refined so that there exists $k \in \{1,\ldots,n\}$ satisfying $\square{\cal C}(g_{z_k})(B_i) \subset \R_+$. We then set $\tilde{\gamma}_i \in \Q$ to the midpoint of the refined box $B_i$ and we add the factor $z_k-g_{z_k}(\tilde{\gamma}_i)$ to the polynomial $\tilde{r}_{k}^{(s)}$. We finally obtain a set of stable univariate polynomials $\tilde{r}_{k}^{(s)} \in \Q[z_k]$, $k=1,\ldots,n$ such that $\sum_{k=1}^{n}{\deg(\tilde{r}_{k}^{(s)})}=d$. 
	
	Let us now introduce the polynomial $\tilde{s}=\prod_{k=1}^{n}{\tilde{r}_{k}^{(s)}}$. By construction $\tilde{s} \in \Q[z_1,\ldots,z_n]$ has rational coefficients and it vanishes on $V(\tilde{I}_r)$, where the polynomial ideal $\tilde{I}_r$ is defined by $\tilde{I}_r=\langle \tilde{f}(t), z_1-g_{z_1}(t),\ldots, z_n-g_{z_n}(t) \rangle$ with $\tilde{f}(t)=\prod_{i=1}^{d}{(t-\tilde{\gamma}_i)} \in \Q[t]$. Hence, according to the classical {\em Nullstellensatz theorem} (\cite{cox1992ideals}), $\tilde{s}$ belongs to the ideal $\tilde{I}_r$ so that there exist polynomials $\tilde{h}_0,\tilde{h}_1,\ldots,\tilde{h}_n \in \Q[t,z_1,\ldots,z_n]$ such that $\tilde{s}=\tilde{h}_0\,\tilde{f}+\sum_{k=1}^n{\tilde{h}_k\,(z_k-g_{z_k})}$. Moreover $\tilde{h}_0$ can be explicitly computed as the quotient of the Euclidean division of $\tilde{s}(g_{z_1}(t),\ldots,g_{z_n}(t))$ by $\tilde{f}(t)$ in $\Q[t]$.
	
	We shall now show that if we refine enough the boxes $B_i$'s isolating the roots $\gamma_i$'s of $f$, then the stable polynomial $s \in I$ that we are seeking for can be obtained from the polynomials $\tilde{s}$ and $\tilde{h}_0$ constructed as explained above. For $\epsilon>0$\footnote{small enough so that the previous process can be applied.}, we denote $\tilde{\gamma}_{i,\epsilon}$, $\tilde{s}_\epsilon$, $\tilde{f}_\epsilon(t)=\prod_{i=1}^d{(t-\tilde{\gamma}_{i,\epsilon})}$, and $\tilde{h}_{i,\epsilon}$ the objects constructed by the previous process where the roots of $f$ are isolated up to precision $\epsilon$ (i.e., $w(B_i) \leq \epsilon$, for all $i \in \{1,\ldots,d\}$). Using the previous notations, the main result of this section can be stated as follows:
	\begin{thm}\label{thm:stabilization} 	
		The polynomial $s=\tilde{s}_\epsilon-\tilde{h}_{0,\epsilon}\,(\tilde{f}_\epsilon-f)$ belongs to the ideal $I_r$. 
		Moreover, there exists $\epsilon>0$ such that the polynomial $s(\sum_{i=1}^{n}{a_i\,z_i},z_1,\ldots,z_n)$ is a stable polynomial.
	\end{thm}
	The proof of Theorem~\ref{thm:stabilization}, given below, requires the following lemma.
	\begin{lem}\label{thm:bound}
		For $0<\epsilon<1$, the polynomial $\tilde{h}_{0,\epsilon}\,(\tilde{f}_\epsilon-f)$ has coefficients bounded by $\epsilon \, \rho$, where $\rho$ is a positive real number that does not depend on $\epsilon$.
	\end{lem}
	\begin{pf}
		Let $f=\sum_{i=0}^{d}{a_i\,t^i}$ and $\tilde{f}_\epsilon=\sum_{i=0}^{d}{b_i\,t^i}$, with $a_d=b_d=1$, denote the expansion of the polynomials $f$ and $\tilde{f}_\epsilon$ on the monomial basis. 
		By the standard \emph{Vieta's formulas}, for all $i \in \{1,\dots,d-1\}$, we have:  
		$$\displaystyle a_i-b_i=\sum_{1\leq k_1 < \cdots < k_i \leq d}{(\gamma_{k_1}\,\cdots\,\gamma_{k_i})-(\tilde{\gamma}_{k_1,\epsilon}\,\cdots\,\tilde{\gamma}_{k_i,\epsilon})}.$$ 
		By assumption, $|\gamma_i-\tilde{\gamma}_{i,\epsilon}| \leq \epsilon$ for all $i \in \{1,\ldots,d\}$, so that 
		{\small $$|a_i-b_i| \leq \left|\sum_{1\leq k_1 < \cdots < k_i \leq d}{(\gamma_{k_1}\,\cdots\,\gamma_{k_i})-(({\gamma}_{k_1}-\epsilon)\,\cdots\,({\gamma}_{k_i}-\epsilon))}\right|.$$}
		Now, we can write $({\gamma}_{k_1}-\epsilon)\,\cdots\,({\gamma}_{k_i}-\epsilon)=\sum_{l=0}^i{\sigma_l\,\epsilon^l}$, where the $\sigma_l$'s denote the \emph{symmetric functions} associated to $\gamma_{k_1},\ldots,\gamma_{k_i}$ and, in particular, $\sigma_0=\gamma_{k_1}\,\ldots\,\gamma_{k_i}$. Consequently, since $\epsilon<1$, we get: 
		{\small $$ |a_i-b_i|\leq\epsilon\,\left|\sum_{1\leq k_1 < \cdots < k_i \leq d}{\sum_{l=1}^i{\sigma_l\,\epsilon^{l-1}}}\right|\leq\epsilon\,\underbrace{\left|\sum_{1\leq k_1 < \cdots < k_i \leq d}{\sum_{l=1}^i{\sigma_l}}\right|}_{\rho_1}.$$}
		On the other hand, the polynomial $\tilde{h}_{0,\epsilon}$ can be computed as the quotient of the Euclidean division of $\tilde{s}_\epsilon(g_{z_1}(t),\ldots,g_{z_n}(t))$ by $\tilde{f}_\epsilon(t)$. Formally, these polynomials can be considered as polynomials in $\Q[\tilde{\gamma}_{1,\epsilon},\ldots,\tilde{\gamma}_{d,\epsilon}][t]$, where $\tilde{\gamma}_{1,\epsilon},\ldots,\tilde{\gamma}_{d,\epsilon}$ are considered as new indeterminates so that their quotient denoted by $h_0$ can be computed independently from $\epsilon$. The coefficients of $h_0$ are thus bounded by a certain positive real number $\delta_1$ that does not depend on $\epsilon$ (use, for instance, \emph{Mignotte's bound} \cite{mignotte1989}). Now, for $i \in \{1,\ldots,d\}$, since $\epsilon<1$, we have $|\tilde{\gamma}_{i,\epsilon}| < |\gamma_i+\epsilon| <|\gamma_i+1|$. Thus if we denote $\delta_2=\max_{i=1,\ldots,d}{|\gamma_i+1|}$, then the coefficients of the evaluation $\tilde{h}_{0,\epsilon} \in \Q[t]$ of $h_0$ for particular values of $\tilde{\gamma}_{1,\epsilon},\ldots,\tilde{\gamma}_{d,\epsilon}$ in $\Q$ are bounded by $\rho_2=\delta_1\,\delta_2^d$.
		Finally, we have proved that the coefficients of $h_{0,\epsilon}(\tilde{f}_\epsilon-f)$ are bounded by $\epsilon \, \rho_1\,\rho_2$, which ends the proof.
	\end{pf}
	We are now in position to give a proof of Theorem~\ref{thm:stabilization}.
	\begin{pf} 
		With the previous notations, for all $\epsilon>0$, we have $s=\tilde{s}_\epsilon-\tilde{h}_{0,\epsilon}\,(\tilde{f}_\epsilon-f)=\sum_{k=1}^n{\tilde{h}_{k,\epsilon}\,(z_k-g_{z_k})} + \tilde{h}_{0,\epsilon}\,f$ so that $s$ vanishes on $V(I_r)$, which implies $s \in I_r$. Let us now prove that we can choose $\epsilon$ so that $s$, viewed as $s(\sum_{i=1}^{n}{a_i\,z_i},z_1,\ldots,z_n)  \in \Q[z_1,\ldots,z_n]$, is stable, i.e., $\forall \lambda \in \Ubar^n$, $|s(\lambda)| >0$. According to Lemma~\ref{thm:bound}, for $\epsilon < 1$, we have $|h_{0,\epsilon}(\lambda)(\tilde{f}_\epsilon(\lambda)-f(\lambda))| \leq \epsilon \, \rho \, \delta$, where $\rho$ (resp., $\delta$) denotes a positive real number (resp., the degree of the polynomial $h_{0,\epsilon}(\tilde{f}_\epsilon-f)$) that does not depend on~$\epsilon$. On the other hand, the polynomial  $\tilde{s}_\epsilon \in \Q[z_1,\ldots,z_n]$ can be written as $\tilde{s}_\epsilon =\prod_{k=1}^{n}{\tilde{r}_{k}^{(s)}} = \prod_{i=1}^d (z_{k_i} - g_{z_{k_i}}(\tilde{\gamma}_{i,\epsilon}))$ for some $k_i \in \{1,\ldots,n\}$ not necessarily distinct. Thus, for all $\lambda \in \Ubar^n$, we have $|\tilde{s}_\epsilon(\lambda)| \geq \prod_{i=1}^d \left| |\lambda_{k_i}| - |g_{z_{k_i}}(\tilde{\gamma}_{i,\epsilon})|\right|$. Now, if we denote by $\displaystyle m=\min_{z \in V(I)}(|z|_\infty-1)$, the minimum distance between the elements of $V(I)$ and $\Ubar^n$, then we have $\left| |\lambda_{k_i}| - |g_{z_{k_i}}(\tilde{\gamma}_{i,\epsilon})|\right| \geq m - \epsilon$, which yields: 
		$$\forall \ \lambda \in \Ubar^n, \ |\tilde{s}_\epsilon(\lambda)| \geq (m-\epsilon)^d.$$ 
		Finally, for sufficiently small $\epsilon$, we have $(m-\epsilon)^d > \epsilon \, \rho \, \delta$ so that:  
		\[
			\begin{aligned}
				\forall \lambda \in \Ubar^n, \ |s(\lambda)| & \geq |\tilde{s}_\epsilon(\lambda)| - |h_{0,\epsilon}(\lambda)(\tilde{f}_\epsilon(\lambda)-f(\lambda))|\\
				& \geq (m-\epsilon)^d - \epsilon \, \rho \, \delta >0,
			\end{aligned}
		\]
		which ends the proof.
	\end{pf}
	The following {\tt StablePolynomial} algorithm summarizes our method for computing a stable polynomial in a zero-dimensional ideal $I$ satisfying $V(I) \cap \Ubar^n= \emptyset$. The routine {\tt IsStable} is used to test if a polynomial $p \in \Q[z_1,\ldots,z_n]$ is stable, i.e., if $V(p) \cap \Ubar^n=\emptyset$ (see \cite{Bouzidi15,BR16}).  
	
	\begin{algorithm}[h]\label{algo:stabilization}
		\caption{{\tt StablePolynomial}}
		{\bf Input:} $I:= \langle p_1,\ldots,p_r \rangle$ be such that  $V(I) \cap \Ubar^n=\emptyset$.
		
		{\bf Output:} $s \in I$ such that $V(s) \cap \Ubar^n=\emptyset$.
		
		{\bf Begin}
		
		$\diamond$ $\{f,g_{z_1},\ldots,g_{z_n}\}:=$ {\tt Univ\_R}($\{p_1,\ldots,p_r\}$);
		
		$\diamond$ $\{B_1,\ldots,B_d\}:=$ {\tt Isolate}($f$);
		
		$\diamond$ $\epsilon:=\min_{i=1,\ldots,d}{w(B_i)}$;		 
		
		{\bf Do}
		
		\hspace{0.3cm} $\diamond$ $[r_1,\ldots,r_n]:=[1,\ldots,1]$ and $\tilde{f} := 1$; 
		
		\hspace{0.3cm} $\diamond$ outside := {\tt False};
		
		\hspace{0.3cm} {\bf For each} $B$ in $\{B_1,\ldots,B_d\}$ {\bf do}
		
		\hspace{0.6cm} {\bf While} (outside={\tt False}) {\bf do}
		
		\hspace{0.9cm} {\bf For} $i$ from $1$ to $n$ {\bf do}
		
		\hspace{1.2cm} {\bf If} $\square{\cal C}(g_{z_k})(B) \subset \R_+$ then
		
		\hspace{1.5cm} $\diamond$ $\gamma:=\text{\tt midpoint}(B)$;
		
		\hspace{1.5cm} $\diamond$ $r_i:=r_i \, (z_i-g_{z_i}(\gamma))$;
		
		\hspace{1.5cm} $\diamond$ outside := {\tt True} and {\bf Break For};
		
		
		\hspace{1.2cm} {\bf End If}
		
		\hspace{0.9cm} {\bf End For}
		
		\hspace{0.9cm} $\diamond$ $\epsilon := \epsilon/2$;
		
		\hspace{0.9cm} $\diamond$ $B:=${\tt Isolate}($f,B,\epsilon$);
		
		\hspace{0.6cm} {\bf End While}
		
		\hspace{0.6cm} $\diamond$ $\tilde{f} := \tilde{f} \, (t-\gamma)$;
		
		\hspace{0.6cm} $\diamond$ outside := {\tt False}; 
		
		\hspace{0.3cm} {\bf End ForEach}
		
		\hspace{0.3cm} $\diamond$ $\tilde{s} := \prod_{i=1}^n{r_i}$;
		
		\hspace{0.3cm} $\diamond$ $\tilde{s}_t := \tilde s$ evaluated at $z_i=g_{z_i}(t)$;
		
		\hspace{0.3cm} $\diamond$ $h_0$ := {\tt quotient}($\tilde{s}_t$,$\tilde{f}$) in $\Q[t]$;
		
		\hspace{0.3cm} $\diamond$ $s := \tilde{s} - h_0\,(\tilde{f}-f)$ evaluated at $t=\sum_{k=1}^n a_k\,z_k$;
		
		{\bf While} ({\tt IsStable}($s$)={\tt False})
		
		$\diamond$ Return $s$.
		
		{\bf End}
		
	\end{algorithm}

	\section{Examples and experiments}\label{sec:experiments}
	
	Let us illustrate the algorithm of Section~\ref{sec:stabilization} on the following simple example:
	$$I=\langle p_1, p_2\rangle, \quad p_1=z_1^2-2\,z_1-2,\; p_2=z_1+z_2-2.$$
	The associated variety $V(I)$ contains two elements, namely, $(1-\sqrt{3}, 1+\sqrt{3})$ and $(1+\sqrt{3}, 1-\sqrt{3})$, so that the stabilizability condition $V(I) \cap \Ubar^2 = \emptyset$ is clearly fulfilled. Yet, $p_1$ and $p_2$ are both unstable polynomials. For $i \in \{1,2\}$, the univariate elimination polynomials (i.e., the resultants of $p_1$ and $p_2$)  $r_{z_i} \in \Q[z_i]$ is given by $r_{z_i} =  z_i^2-2\,z_i-2$. The polynomial $z_i^2-2\,z_i-2$ being irreducible in $\Q[z_i]$, this makes the approach of \cite{xu1994} impracticable. Let us apply the algorithm of Section~\ref{sec:stabilization} for computing a stable polynomial $s \in I$. We start by computing a univariate representation of $V(I)$. We get:
	$$ f(t) := t^2-2\,t-2  = 0, \quad  z_1 = t, \quad  z_2 = 2-t.$$
	The roots of $f(t)$ are given by $\gamma_1 \approx -0.73$ and $\gamma_2 \approx 2.73$ and choosing the precision  $\epsilon = \frac 1 2$, we get the approximate roots (in $\Q$) $\tilde{\gamma}_1=-\frac{1}{2}$ and $\tilde{\gamma}_2=3$. Consequently, the algorithm of Section~\ref{sec:stabilization} yields
	$$ \tilde{f}(t)= \left(t+\frac{1}{2}\right)\,(t-3), \quad \tilde{s}(z_1,z_2) = (z_1-3)\,\left(z_2-\frac{5}{2}\right),$$
	which then leads to:
	$$h_0(t)= -1, \quad (\tilde{f}-f)(t) = -\frac{1}{2} \, t  + \frac{1}{2}.$$
	Finally, after substituting $t=z_1$ in $\tilde{f}-f$, we get:
	$$s(z_1,z_2)  = z_1\,z_2-3\,z_1-3\,z_2+8.$$
	We can then check that this polynomial is stable so that we are done (see \cite{Bouzidi15}). As a byproduct, we also obtain the corresponding cofactors, that is, $$s=-p_1+(z_1-3)\,p_2 \in I.$$
	
\begin{figure}
		\includegraphics[width=\columnwidth]{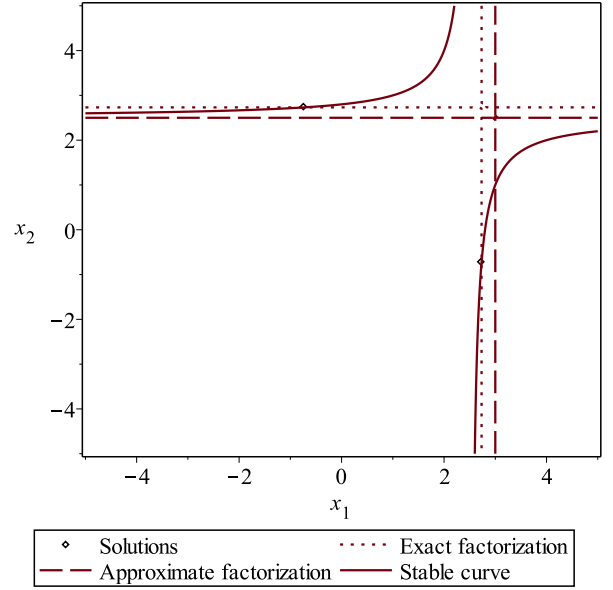}
		\caption{Stabilizing polynomial for the variety $V(I)$ corresponding to $I=\langle z_1^2-2\,z_1-2, \, z_2^2-2\,z_2-2\rangle$}\label{fig:example}
	\end{figure}

	Figure~\ref{fig:example} shows the stable polynomial obtained by an exact factorization of $r_{z_1}$ and $r_{z_2}$ in dots, the approximate factorization $\tilde{s}$ used in our algorithm in dash, and finally the stable polynomial $s$ obtained after adding the correcting term represented by the solid curve.
	
	We have implemented two routines {\tt IsStabilizable} and {\tt StablePolynomial}, which correspond respectively to the algorithms given in Sections~\ref{sec:stabilizability} and \ref{sec:stabilization}, in the computer algebra system {\sc Maple}. These routines use the procedure {\tt resultant} for computing the resultant of two polynomials (required in the computation of $l_k$ in Algorithm~\ref{algo:stabilizability}), the procedure {\tt RationalUnivariateRepresentation} of the {\sc Maple} package {\tt Groebner} for computing the \emph{Rational Univariate Representation}\footnote{This rational representation, which outputs expressions for the coordinates that are rational fractions, is post-processed in order to get polynomial expression for the coordinates as defined in Definition~\ref{def:UR}.}
	of a zero-dimensional ideals, and the procedure {\tt fsolve} for computing the complex roots of a univariate polynomial \footnote{The detailed code along with the used testsuite can be found in \cite{online}}.  
	
	In Table~\ref{table1} below, we report the running times (in seconds) of {\tt IsStabilizable} and {\tt StablePolynomial} applied to systems of randomly chosen polynomials in two or three variables with integer coefficients chosen uniformly at random between $-100$ and $100$. In order to get zero-dimensional systems, we choose as many polynomials as number of variables.
	Moreover, we use the change of variables $z_i = Z_i/10$, $i=1,2$ or $i=1,2,3$ to increase the probability of the roots to be outside the unit polydisc. The experiments have been conducted on 2.10 GHz Core(TM) Intel i7-4600U with 4MB of L3 cache
	with {\sc Maple 2015} under windows platform.
	
	\begin{table}
		\begin{center}
			\begin{tabular}{|c|c||c|c|}
				\hline
				\multicolumn{2}{|c||}{Data} & \multicolumn{2}{c|}{Running time}  \\[5pt]
				\hline
				nbvar &  $\#V(I)$  & {\tt IsStabilizable} & {\tt StablePolynomial}  \\[5pt]
				\hline
				\multirow{4}{1cm}{2}  & 9 & 0.09 &  0.11 \\
				\cline{2-4}
				&  25 & 1.23 & 0.50\\
				\cline{2-4}
				& 64 & 38.10 & 7.84\\
				\cline{2-4}
				& 100 & 244.91 & 49.49 \\
				\hline
				\multirow{4}{1cm}{3} & 8 & 0.13 &  0.11\\
				\cline{2-4}
				& 27 & 4.39 & 0.87\\
				\cline{2-4}
				& 36 & 11.83 & 1.98 \\
				\cline{2-4}
				& 48 & 33.92 & 5.32\\
				\cline{2-4}
				& 64 & 118.28 & 24.09\\
				\hline
			\end{tabular}
		\end{center}
		\medskip
		\caption{CPU times in seconds of {\tt IsStabilizable} and {\tt
				StablePolynomial} runned on sets of random polynomials in $2$ or $3$
			variables with integer coefficients.}\label{table1}
	\end{table}

	\begin{remark}
		From Table~\ref{table1}, one can notice that, in general, the running times of {\tt IsStabilizable} are higher than those of {\tt StablePolynomial}. This is, most likely, due to the additional cost induced by the computation of~$l_k$ in Algorithm~\ref{algo:stabilizability}, which requires the computation of  elimination polynomials for each variable (resultants). 
		\end{remark}

	
	\bibliography{ifacconf}  
	%
	%
	
\end{document}